\documentclass{amsart}[15]
\usepackage{latexsym,amssymb,graphicx,amsmath,amscd}

\usepackage[all,cmtip]{xy}

\usepackage[mac]{inputenc}

\usepackage{amssymb}

\newcommand{\CC}{\Bbb C}
\newcommand{\PP}{\Bbb P}

\newtheorem{defi}{Definition}[section]
\newtheorem{pro}[defi]{Proposition}

{\newtheorem{lem}[defi]{Lemma}

\newtheorem{theo}[defi]{Theorem}
\newtheorem{coro}[defi]{Corollary}
\newtheorem{exa}[defi]{Example}

\newtheorem{rem}[defi]{Remark}
\newtheorem{rems}[defi]{Remarks}
\newtheorem{conj}[defi]{Conjecture}

\title [On the smooth locus of the aligned Hilbert Scheme ]{ \hspace{0,4cm} On the smooth locus of aligned Hilbert Schemes  \hspace{3cm}
 \hspace{1cm}The $k$-secant lemma  and the general projection theorem  }

\author{Laurent Gruson }
\address{Laboratoire de Mathématques de Versailles (UMR 8100 CNRS)\\
45, avenue des États-Unis\\
F-78035 Versailles cedex
France}
\email{laurent.gruson@math.uvsq.fr}

\author{Christian Peskine}
\address{Institut de Mathématiques de Jussieu (UMR 7586 CNRS)\\  Université Pierre et Marie Curie  \\Case 247\\4, place Jussieu\\75252 Paris Cedex 05\\
FRANCE }
\email{peskine@math.jussieu.fr}
\begin{document}

\date{October 10, 2010}

\begin{abstract} Let  $X$ be a smooth, connected, dimension $n$, quasi-projective variety embedded in  ${\PP}^{N}$.
Consider integers  $\{k_1,...,k_r\}$, with $k_i>0$, and the Hilbert Scheme  $H_{\{k_1,...,k_r\}}(X) $   of aligned,  finite, degree $\sum k_i$, subschemes of  $X $, with multiplicities $k_i$ at points $x_i$ (possibly coinciding). The expected dimension of $H_{\{k_1,...,k_r\}}(X) $ is
 $2N-2+r-(\sum k_i)(N-n)$. We study the locus of points where $H_{\{k_1,...,k_r\}}(X) $  is not smooth of expected dimension and
we prove that the lines carrying this  locus do not fill up ${\PP}_{N}$.

 \end{abstract}

\maketitle

 
 
 
 \section{Introduction}

 \vspace{0,5cm}
 
Let $C \subset {\PP}^3({\CC})$  be a smooth curve in the projective complex space . A general projection
$p : C \rightarrow C_1  \subset  {\PP}^2({\CC})$  has only ordinary double points as singularities. This statement, known as the $3$-secant lemma, is composed of three assertions:

1) the tangents to $C$ do not fill up the space,

2) the tacnode, or stationary, or ramified $2$-secant lines to $C$ do not fill up the space, 

3)  the  $3$-secants to $C$ do not fill up the space.

The proof is classical and easy to explain . We note that 1) is obvious (counting dimensions). If 2) were not true, two tangents would always intersect.  Consequently, if $C$ is not a plane curve,  all tangents would pass through a point and C would be everywhere ramified over its projection from this point. As for 3), it reduces to 2). Indeed, if every $2$-secant to $C$ is a $3$-secant to $C$,  it is not difficult to check that two tangents always intersect.

It is well known that the  double locus $C_2$ of the projection $C_1$  has a natural structure of smooth variety  whose ideal in $C_1$ is the conductor. The tangent space to $C_2$ is implicitly described in the  $3$-secant lemma. Consider  $z\in C_2$ and the points $x_1,x_2\in p^{-1}(z)$, the tangent space to $C_2$ at $z$ is the intersection of the projections of the tangent spaces (lines) to $C$ at $x_1$ and $x_2$ (they intersect transversally).

Before  discussing  possible generalizations of this  result to  higher dimensions, let us agree that in this paper a line $L \subset {\PP}^{N}({\CC})$ is a $k$-secant to a smooth quasi-projective variety $Z \subset {\PP}^{N}({\CC})$ if the scheme $L\cap Z$ is finite of degree  
$\geq k$.

The $3$-secant lemma was first generalized  by Z. Ran (\cite{8}) as follows: the $n+2$-secants to a smooth, dimension $n$,  projective variety 
$X \subset {\PP}^{N}({\CC})$ fill up a variety of dimension at most $n+1$.

Recently R. Beheshti and D. Eisenbud improved significantly Ran's  lemma (see \cite{3}, Theorem 1.5.):  they prove that for 
$k >[n/s] + 1$ (where $[n/s]$ is the integral part of $n/s$), the $k$-secant lines to a smooth, dimension $n$,  projective variety $X \subset {\PP}^{N}({\CC})$ fill up a variety of dimension at most $n+s$.

Note that if $c=N-n$ and if $k \leq n/(c-1)+1$, the $k$-secants to $X$ may fill up the space. In this case, the projection of $X$ from a general point may have points of order $\geq k$ and we need information such as  dimension, smoothness, description of tangent spaces  concerning the geometric nature  of their  locus.  For example, using the result of Beheshti/Eisenbud, it is clear that  the projection (from a general point) of a smooth, dimension $6$, variety $X \subset {\PP}^{9}({\CC})$ has no points of order $5$. We would like to know  the dimension and  the singular locus of the loci of points of order $k$, for $2\leq k \leq 4$ for this projection.

\vspace{0,3cm}

Here is our first result.

\begin{theo} (General Projection Theorem) Let  $X \subset {\PP}^{N}$ be a smooth variety of dimension $n$  and codimension $c=N-n$,   and let $\pi : X \rightarrow X_1 \subset {\PP}^{N-1}$ be a projection from a general point of ${\PP}^{N}$.

1) For $k=k_1+...+k_r$, with $k_i>0$, let $X_{\{k_1,...,k_r\}} \subset X_1$ be the subscheme formed by points $x\in X_1$ such that $\pi^{-1}(x)$ contains $r$  points $\{x_1,...,x_r\}$  (possibly coinciding) with  multiplicity $\geq k_i$ in $x_i$. Then

1) The scheme $ X_{\{k_1,...,k_r\}}$ has pure dimension $N-1-\sum_{i=1}^r  (k_ic-1) =N-1+r-kc$ (the empty set has all dimensions).
  
 2) The singular locus of $X_{\{k_1,...,k_r\}}$ is $X_{\{k_1,...,k_r,1\}}$.

3) The normalization $\tilde{X}_{\{k_1,...,k_r\}}$ of $X_{\{k_1,...,k_r\}}$ is smooth.

\label{genproj}

\end{theo}

\begin{rems}

1) Be careful, if $x_i$ and $x_j$ do coincide, the multiplicity of $\pi^{-1}(x)$  at the point  $x_i=x_j$ has to be $\geq (k_i+k_j)$.

2) Please note the following special case of this theorem.

When $k_i=1$ for all $i$, the scheme  $X_k=X_{\{1,...,1\}} \subset X_1$ formed by points of multiplicity $\geq k $ of $X_1$ has dimension  $N-1-k(c-1)$.  The singular locus of $X_k$  is $X_{k+1}$ and the normalization $\tilde{X_k}$ of $X_k$  is smooth.
\end{rems}

\vspace{0,2cm}

It is perhaps worthwhile to emphasize here that the main difficulties to generalize the $3$-secant lemma to any dimension appear when the tangent spaces to $X$ fill up the ambient space.

Indeed, suppose they don't. Then the fiber of a point $x \in X_k$ is reduced. By a simple computation in the Grassmann variety $G(1,N)$, one sees that $X_k$ is smooth of dimension $N-1-k(c-1)$ at a point $x$ if and only if the fiber of $x$ has degree $k$ and  the projections of the tangent spaces to $X$ at the $k$ distinct points of the fiber  are in relatively general position, in which case the tangent space to $X_k$ at $x$ is their intersection. Assuming they are not, there is a corresponding special configuration of  linear spaces contained in the Segre ${\PP}_1\times {\PP}_{N-2}$ whose projective cone is the intersection of $G(1,N)$ with its tangent space at the point (line) corresponding to $x$. Imitating the proof of the $3$-secant Lemma, one can prove that if these  linear spaces (in $k$ distinct ${\PP}^{N-2}$ of the Segre) are not in relative general position, the projections of any $k-1$ among them are not  in relatively general position.  The conclusion  of the proof  in this case goes through an easy  analysis of  aligned Hilbert schemes (see section $2$ below).

If the tangent spaces to $X$ fill up the ambient space, one can consider the open subset $\hat{X}_k \subset X_k$   formed by points $x \in X_k$ whose fiber is reduced of degree $k$. The same argument proves that $X_k$ is smooth of dimension $N-1-k(c-1)$ at the points of $\hat{X}_k$.

There is no intuitive geometric description of the tangent space to $X_k$ in a point corresponding to a tangent line to $X$ (except for $k=2$). That is why to prove our main results we have to describe algebraically (and with brutal force) the local equations and the local cotangential equations of $X_k$. 

To our knowledge, it was not known that $\hat{X}_k $ is a dense open set in $X_k$. This is implicit in our Theorem \ref{genproj}, in particular
$$\hat{X}_k =\emptyset \Rightarrow X_k=\emptyset. $$

\vspace{0,3cm}

Our  General Projection Theorem  is not detailed enough. We should be more precise about the scheme structure of the closed algebraic set  $X_{\{k_1,...,k_r\}}$. The best way to do this, and more generally to clarify our point of view,  is  to state and  prove our result in the language of Hilbert Schemes of aligned points.  We know that such Hilbert Schemes are well defined and equipped with an obvious map to the Grassmann variety of lines. The normalization  $\tilde{X}_{\{k_1,...,k_r\}}$ is the inverse image of $X_{\{k_1,...,k_r\}}$ in the corresponding  Hilbert Scheme.
 This is explained and described in the following theorem, of which the previous one is clearly a consequence. 

From now on, $G(1,N)$ is the Grassmann variety of lines in ${\PP}^{N}$  and we denote  by $ {\mathcal  I} \subset G \times {\PP}^{N}$ the incidence variety point/line. We recall that ${\mathcal  I} $ is a projective line bundle over $G$ on one hand, and  a (${\PP}^{N-1}$)-bundle over ${\PP}^{N}$ on the other hand.

\begin{theo} (Aligned Hilbert Scheme Theorem)   Let  $X$ be a smooth, connected, dimension $n$, quasi-projective variety embedded in  
${\PP}^{N}$, with $N=n+c$.

For $k=k_1+...+k_r$, with $k_i>0$, let $H_{\{k_1,...,k_r\}}(X) $ be the  Hilbert scheme  of aligned,  finite, degree $k$ subschemes of $X$,
with multiplicities $k_i$ in points $x_i$ (possibly coinciding).  Consider  the natural projective line bundle $H_{\{k_1,...,k_r\}}(X) \times_{G}  {\mathcal  I} $   over
 $H_{\{k_1,...,k_r\}}(X)$ and the projection 
 $$\theta_{\{k_1,...,k_r\}} : H_{\{k_1,...,k_r\}}(X) \times_{G}  {\mathcal  I} \rightarrow {\mathcal  I}  \rightarrow  {\PP}^{N}.$$
 Then the general fiber of  $\theta_{\{k_1,...,k_r\}}  $ is smooth of pure dimension $N-1+r-kc$.
 
 \label{hilbsch}
\end{theo}

As in the case of the general projection theorem, the following remarks are important.

\begin{rems}
1) When $x_i$ and $x_j$  coincide, the multiplicity of a point  $h \in H_{\{k_1,...,k_r\}}(X)$  at   $x_i=x_j$ has to be $(k_i+k_j)$.

2) Please note the following special case of this theorem.

When $k_i=1$ for all $i$, if we denote 
$H_k(X)=H_{\{1,...,1\}}(X)$  
the Hilbert  scheme of aligned,  finite, degree $k$ subschemes of $X$ and 
$\theta_k : H_k(X) \times_{G}  {\mathcal  I} \rightarrow {\mathcal  I} \rightarrow  {\PP}_{N}$, then the general fiber of $\theta_k $ is smooth of dimension $N-1-k(c-1)$.

\end{rems}

\vspace{0,2cm}

As a  special case of  Theorem  \ref{hilbsch}, for $r=1$ and any $k$,   we recover a well known result  of Mather (see \cite{7}):  ``higher polar varieties"  of a general point with respect to a smooth variety $X$ cut in $X$ a smooth variety of expected dimension.

\begin{coro} (Mather) Let  $HB_k(X) = H_{\{k\}}(X) \subset  {\mathcal  I}$ 
be  the Hilbert-Boardmann locus of all $(L,x) \in   {\mathcal  I}$ such that $ L\cap X$ has multiplicity
at least $k$ at the point $x$.  Consider  the natural projective line bundle $HB_{k}(X) \times_{G}  {\mathcal  I} $ and the projection
 $$\theta_{\{k\}} : HB_{k}(X) \times_{G}  {\mathcal  I} \rightarrow {\mathcal  I}  \rightarrow  {\PP}^{n+c}.$$
 Then the general fiber of  $\theta_{\{k\}}  $ is smooth of pure dimension $n-(k-1)c=N-kc$. 
\end{coro}

We can as well note here that the result of R. Beheshti and D. Eisenbud is recovered as a direct consequence of Theorem  \ref{hilbsch}.
Indeed, assume  $k >[n/s] + 1$ (where $[n/s]$ is the integral part of $n/s$).  We need to show that the $k$-secant lines  to $X$ fill up a variety of dimension at most $n+s$.

For $s\geq c$, there is nothing to prove. If $s=c-1$, this is a special case of our theorem. Assume $s\leq c-2$ and let $L$ be a $k$-secant line of $X$.
Consider a projection $X \rightarrow {\PP}_{n+s+1}$ whose double locus  avoids  the finite scheme $L\cap X$. By  Theorem  \ref{hilbsch}, the $k$-secants to the smooth locus of the image of this projection fill at most a hypersurface in ${\PP}_{n+s+1}$,  hence the $k$-secant lines of $X$  near $L$ fill a variety of dimension at most $n+s$ in ${\PP}^N$.

\qed

\vspace{0,3cm}

Our  Aligned Hilbert Scheme Theorem  is an easy consequence of  the  Aligned Ordered Hilbert Scheme Theorem. 

The ordered Hilbert schemes $OH_{(k_1,...,k_r)}(X)$ parametrizes finite aligned subschemes $Z \subset X$ supported in an ordered set  of points $(x_1,...,x_r) \in X^r$ (not necessarily distinct)   and with  ordered multiplicities $ k_i$ at $x_i$ 
 (note once again that if a point is redundant, for example if $x=x_{i_1}=...=x_{i_s}$, then $Z\subset L\cap X$ must have multiplicity $ k_{i_1}+...+k_{i_s}$ at $x$).

Since $OH_{(k_1,...,k_r)}(X)$ is  finite and flat over $H_{\{k_1,...,k_r\}}(X) $, it is clear that if $OH_{(k_1,...,k_r)}(X)$ is smooth, then so is $H_{\{k_1,...,k_r\}}(X) $ (the converse is not true). Theorem  \ref{hilbsch} is then a straightforward corollary of  the following stronger result.

\begin{theo} (Aligned Ordered Hilbert Scheme Theorem)   Let  $X$ be a smooth, connected, dimension $n$ quasi-projective variety embedded in  
${\PP}^{N}$, with $N=n+c$.

For $k=k_1+...+k_r$, with $k_i>0$, let $OH_{(k_1,...,k_r)}(X) $ be the  ordered Hilbert scheme  of aligned,  finite, degree $k$ subschemes of $X$,
with (ordered) multiplicities $k_i$ at the ordered  points $x_i$ (possibly coinciding).  Consider  the natural projective line bundle 
$H_{(k_1,...,k_r)}(X) \times_{G}  {\mathcal  I} $   over
 $H_{(k_1,...,k_r)}(X)$ and the projection 
 $$\theta_{(k_1,...,k_r)} : OH_{(k_1,...,k_r)}(X) \times_{G}  {\mathcal  I} \rightarrow {\mathcal  I}  \rightarrow  {\PP}{N}.$$
 
 The general fiber of  $\theta_{(k_1,...,k_r)}  $ is smooth of dimension $N-1+r-kc$.
 
 \label{ohilbsch}
\end{theo}

\vspace{0,4cm}

The three following sections are devoted to the proof of this theorem. 
 \vspace{0,2cm}

In section $2$ we consider a closed subvariety $\varGamma$ of an affine line ${\Bbb A}_{ \mathop{\rm Spec}  R}^1=\mathop{\rm Spec}  (R[z])$ over an affine smooth variety, i.e.   $R$ is a finitely generated regular ${\CC}$-algebra. To the morphism

$$\phi : \varGamma= \mathop{\rm Spec}  (R[z]/J) \rightarrow  \mathop{\rm Spec}  R$$
are associated ordered aligned  Hilbert Schemes, that we denote
$OH_{(k_1,...,k_r)}(\varGamma)$ or $OH_{(k_1,...,k_r)}(\phi)$, parametrizing  subschemes of the fibers of $\phi $ with ordered multiplicities $k_i$ in  ordered sets of points of the fibers.
Such Hilbert Schemes are equipped with obvious set maps 
$$OH_{(k_1,...,k_r)}(\phi)=OH_{(k_1,...,k_r)}(\varGamma) \rightarrow 
\mathop{\rm Spec}(R[z_1,...,z_r]).$$
In Proposition \ref{defhilb} we recall   that these maps are embeddings which we describe. We state and prove two general  technical lemmas that we use repeatedly  in the sequel. In particular, Lemma \ref{hilbhyper} describes the local equations and the local cotangential equations of the embedding $OH_{(k_1,...,k_r)}(\phi)\subset \mathop{\rm Spec} (R[z_1,...,z_r])$. This is elementary calculus.
 \vspace{0,2cm}

In section $3$, we focus on the case where the base affine variety is an open affine subvariety is a Grassmann variety. More precisely we interpret the last lemma of section $2$ in two special cases.
On the one hand, when $\mathop{\rm Spec} R$ is an open set of  $G(1,N)$ (the Grassmann  variety of lines in ${\PP}^N$), and on the other hand when $\mathop{\rm Spec} R$
is an open set of ${\PP}^{N-1}(\beta)$ (the  Grassmann  variety of lines through a point $\beta \in {\PP}^N$).
 \vspace{0,2cm}

The proof of  Theorem \ref{ohilbsch} is  presented in section $4$. We follow an induction principle inspired by the classical proof of the $3$-secant lemma. The main difficulty stems from the fact that there is no  natural  geometric description of the tangent space to the Hilbert-Boardmann Scheme $HB_k(X)$ in a general point. We overcome this difficulty by exploiting the ``algebraic Segre nature" of a tangent space to the Grassmann variety of lines. This Segre structure is described with all necessary precautions in section $3$.
\vspace{0,2cm}
 
The last section is dedicated to examples, questions and conjectures.
\vspace{0,2cm}

As a conclusion to this introduction, we wish to thank the referees for their constructive remarks and critics.

\section{The local  Aligned Ordered Hilbert Scheme.} 

In this section $R$ is  a regular  finitely generated ${\CC}$-algebra. We consider an affine line $ \mathop{\rm Spec}R[z]$ over the affine smooth variety $\mathop{\rm Spec}  R$,
 a closed  subscheme  $\varGamma=\mathop{\rm Spec}  (R[z]/J) $  of this affine line and the morphism

$$\phi : \varGamma=\mathop{\rm Spec}  (R[z]/J) \rightarrow  \mathop{\rm Spec}  R.$$

The aligned ordered   Hibert Scheme $OH_{(k_1,...,k_r)}(\phi)$ parametrizes subschemes 
$$Z \subset \mathop{\rm Spec} {\CC}[z]=\phi^{-1}(x), \quad x \in \mathop{\rm Spec}  R$$
with support in an ordered  set of points $a_1,...,a_r \in \mathop{\rm Spec}  {\CC}[z]$,
and with length (multiplicity) $k_i$ at $a_i$.  Two points $a_i$ and $a_j$ may coincide, as long as 
the finite scheme $Z$ has length multiplicity $ k_i+k_j$ in $a_i=a_j$.

When $\varGamma$ is a hypersurface, i.e. when $\varGamma=\mathop{\rm Spec}  (R[z]/(g))$, where $g=g(z)$ is a  polynomial, we often write $ OH_{(k_1,...,k_r)}(g)$ instead of 
$OH_{(k_1,...,k_r)}(\phi) $. 
\vspace{0,2cm}

The following proposition is  well known to anyone familiar with aligned Hilbert schemes (see for example \cite{4} or \cite{5}). For  a reader  who is not,  the best is to admit 1), which by the way explains  why we prefer the ordered Hilbert scheme to the nonordered  one. 

\begin{pro} 

1) We have

$$OH_{(k_1,...,k_r)}(g)=\mathop{\rm Spec} (R[z_1,...,z_r]/(h_0,...,h_{k-1})), \   k=\sum k_i , $$
where the polynomials $h_l \in R[z_1,...,z_r] $ are defined for  $l=0,...,k-1$ by 

$$g(z)\equiv \sum_0^{k-1}  h_l(z_1,...,z_r)z^l \  \mathop{\rm mod} ( \prod_1^r  (z-z_i)^{k_i}).$$

2)   If $\varGamma=\mathop{\rm Spec} (R[z]/(g_1,...,g_c))$, then 

$$OH_{(k_1,...,k_r)}(\phi) = \cap_1^c   H_{k_1,...,k_r}(g_t).$$

3) If $R'$ is a finitely generated regular $R$-algebra and if  $\phi'=\phi \otimes_R R' :  \mathop{\rm Spec} \ A\otimes_R R' \rightarrow  \mathop{\rm Spec} \ R'$, then 
$$OH_{(k_1,...,k_r)}(\phi')=OH_{(k_1,...,k_r)}(\phi)\times_{\mathop{\rm Spec} R} (\mathop{\rm Spec}  R').$$

\label{defhilb}
\end{pro}

\vspace{0,3cm}

We observe (with pleasure) that $OH_{1}(g) \simeq \mathop{\rm Spec} (R[z]/(g(z))$ and more generally  
$OH_{1}(\phi) \simeq  \varGamma $.
\vspace{0,2cm}

From now on, we  shall pay a particular attention to the case  when  $\varGamma = \mathop{\rm Spec} \ (R[z]/J)$ is a smooth complete intersection in 
$ \mathop{\rm Spec} \ R[z]$.

We begin with obvious remarks.

\begin{rems}

1) The expected dimension of $OH_{(k_1,...,k_r)}(g)$ is $ \mathop{\rm dim} R + r -k = \mathrm{dim} \ R +  \sum_1^r (1-k_i)$.

2)  When $R[z]/J$ is a complete intersection of codimension $c$  in $R[z]$,  the expected dimension of $OH_{(k_1,...,k_r)}(\phi) $ is $\mathop{\rm dim} \ R + r -ck$.

\end{rems}

We note  (with great pleasure once  again) that since $OH_1(\phi) \simeq  \varGamma$ it is clear that when 
$\varGamma $ is smooth, so is the ordered Hilbert Scheme  $OH_1(\phi)$. This obvious  remark will be the starting point of the proof by induction of Theorem
 \ref{ohilbsch}.
\vspace{0,2cm}

The following result will prove to be an important technical tool in the proof of our main theorem. To be more precise, it will allow us, when necessary,
to work with points $(x,a_1,...,a_r) \in OH_{(k_1,...,k_r)}(\phi)\subset \mathrm{Spec} \ R[z_1,...,z_r] $, with $x\in  \mathop{\rm Spec} R$ and $a_i\in {\CC}$,  such that $a_i \neq a_j$ for $i \neq j$.

\begin{lem}  Assume $R[z]/J$ is a complete intersection of codimension $c$ in $R[z]$.

Consider a point $(x,a_1,...,a_r) \in OH_{(k_1,...,k_r)}(\phi) $ (with $x\in \mathop{\rm Spec}  R$ and  $a_1,...,a_r \in \mathop{\rm Spec} \  {\CC}[z]=\phi^{-1}(x)$). Assume that there exist  $1\leq s < t \leq r$ 
such that $a_s=a_t$, in other words that
$$(x,a_1,...,a_s,...,a_{t-1},a_{t+1},...,a_r) \in OH_{(k'_1,...,k'_{r-1})}(\phi),$$
with 
  $k_i=k'_i $ for $ i<s $ and $ s<i<t $,   $k'_s=k_s+k_t $  and  $ k'_{i-1}=k_{i}$  for  $ i>t.$

Then   $OH_{(k_1,...,k_r)}(\phi) $ is smooth of expected dimension at $(x,a_1,...,a_r)$ if and only if  $OH_{(k'_1,...,k'_{r-1})}(\phi) $ is smooth of expected dimension at $(x,a_1,...,a_s,...,a_{t-1},a_{t+1},...,a_r)$.

\label{lpointsdistinct}

\end{lem}

\proof

\vspace{0,3cm} 

For the sake of simplicity,  we assume 
$$R[z]/J=R[z]/(g), \quad  s=1, \quad   r=k=t=2, \quad  k_1=k_2=1, \quad a_1=a_2=0.$$
 
Put $g(z)= \sum_{i\geq 0} \alpha_i z^{d-i}$. Let  ${\mathcal M}$ be the maximal ideal of $R$ corresponding to the point $x\in \mathrm{Spec} \ R$. 

Since $0$ is a point of multiplicity $\geq 2$ in the fiber of $x$, we have $ \alpha_{d},  \alpha_{d-1} \in {\mathcal M}$.

$$ g(z) \equiv  (\alpha_{d-1}+ \alpha_{d-2}(z_1+z_2))z+ \alpha_d  \  \ \mathop{\rm mod } ( (z-z_1)(z-z_2)+({\mathcal M}, z_1, z_2)^2)R[z_1,z_2][z], $$

and 

$$ g(z) \equiv ( \alpha_{d-1}+ 2\alpha_{d-2}z_1)z+ \alpha_d  \ \  \mathop{\rm mod } ((z-z_1)^2+ ({\mathcal M}, z_1)^2)R[z_1][z].$$

\vspace{0,2cm} 

Two  cases occur (depending of the multiplicity of the root $0$ of the image of $g$ in $(R/{\mathcal M})[z]={\CC}[z]$).
\vspace{0,2cm} 

1) If the  multipicity is precisely $2$, i.e. if $ \alpha_{d-2} \notin {\mathcal M}$, then 

$$OH_{1,1}(g) \ \mbox{ is smooth of expected dimension } \mathrm{dim} \ R  \Leftrightarrow  \alpha_d \notin {\mathcal M}^2$$
$$  \Leftrightarrow
OH_{2}(g) \ \mbox{  is smooth of expected dimension } \mathrm{dim} \ R-1.$$ 

2)  If the  multipicity is  $>2$, i.e if $ \alpha_{d-2} \in {\mathcal M}$, then 

$$OH_{1,1}(g) \ \mbox{ is smooth of expected dimension } \mathrm{dim} \ R$$
$$  \Leftrightarrow  \alpha_{d-1} \mbox{ and }  \alpha_d \mbox{ are transverse in }
{\mathcal M}/ {\mathcal M}^2$$
$$  \Leftrightarrow
OH_{2}(g) \ \mbox{  is smooth of expected dimension } \ \mathrm{dim} \ R-1.$$

\qed

\vspace{0,3cm}

Our next result  describes explicitly the local equations and the local cotangential equations of the ordered Hilbert Scheme at a point 
$$ (x,a_1,...,a_r)  \in OH_{(k_1,...,k_r)}(g) \subset \mathrm{Spec} \ R[z_1,...,z_r],$$
 with $x \in \mathop{\rm Spec} R $ and  $ a_i \in {\CC}$.  It will be used more than once (and without thinking twice).

From now on, we denote by $g^{(s)}(z)$ the derivative of order $s$ of the function $g(z)$ for $s\geq 0$.  The convention  $g^{(-1)}(z)=0$  will prove to be useful later on.

\begin{lem}  Assume that $\varGamma=\mathop{\rm Spec}(R[z]/(g))$ is a hypersurface  and consider a point
$(x,a_1,...,a_r) \in OH_{(k_1,...,k_r)}(g)$ supported in the fiber $\phi^{-1}(x)$, with $x \in \mathop{\rm Spec}  R$  and $a_i \in {\CC} $,  such that $a_i \neq a_j$ for $i\neq j$.

1) The local equations of $OH_{(k_1,...,k_r)}(g) \subset  \mathrm{Spec} \ R[z_1,...,z_r] $, at $(x,a_1,...,a_r) $, are 

$$g^{(s)}(z_i)/s!  ,   \ \  i=1,...,r ,  \ \ 0 \leq s < k_i.$$

2) If ${\mathcal M}$ is the maximal ideal of $R$ corresponding to the point $x \in {Spec} \ R$, the local cotangential equations of $OH_{(k_1,...,k_r)}(g) \subset  \mathop{\rm Spec} R[z_1,...,z_r] $ at $(x,a_1,...,a_r) $  in the cotangent space  
$$({\mathcal M},(z_1-a_1),...,(z_r-a_r))/({\mathcal M},(z_1-a_1),...,(z_r-a_r))^2$$
of $\mathop{\rm Spec} R[z_1,...,z_r] $ in the point $(x,a_1,...,a_r)$ ,  are the classes  of  the $r(\sum k_i)$ elements

$$g^{(s)}(a_i)  \  \  0 \leq s < k_i-1 ,   \ \ \           
 g^{(k_i-1)}(a_i)+(z_i-a_i)g^{(k_i)}(a_i)$$
for $  i=1,...,r$.

\label{hilbhyper}

\end{lem}

\proof

1) is an obvious consequence of the Taylor expansions $g(z) = \sum_{s>0}  (g^{(s)}(z_i)/s!)(z-z_i)^s$.

\vspace{0,3cm}

2) is easily deduced from the relation $g^{(s)}(z_i) \equiv g^{(s)}(a_i)+(z_i-a_i)g^{(s+1)}(a_i) \ \mbox{ mod } (z_i-a_i)^2$.

\qed

\vspace{0,3cm}

The following  remarks (using the same notations as in the lemma) are important.

\begin{rems}

1) For $s<k_i-1$,   the classes  
 $$cl(g^{(s)}(a_i)) \in  ({\mathcal M},(z_1-a_1),...,(z_r-a_r))/({\mathcal M},(z_1-a_1),...,(z_r-a_r))^2$$
are in the  vector subspace  ${\mathcal M}/{\mathcal M}^2$

\vspace{0,3cm}

2) If $g^{k_i}(a_i) \in {\mathcal M}$, i.e. if $cl(g(z)) \in (R/{\mathcal M})[z]$ has multiplicity $>k_i$ at $a_i$, then 
$$cl(g^{(k_i-1)}(a_i)+(z_i-a_i)g^{(k_i)}(a_i))=cl(g^{(k_i-1)}(a_i)) \in {\mathcal M}/{\mathcal M}^2.$$

\vspace{0,3cm}

3) If $g^{(k_i)}(a_i) \notin {\mathcal M}$, i.e.  the order of $g(z)$ at the point $a_i$ of the special fiber is precisely $k_{i}$, then
$$cl(g^{(k_i-1)}(a_i)+(z_i-a_i)g^{(k_i)}(a_i)) \notin {\mathcal M}/{\mathcal M}^2.$$
\end{rems}

  \vspace{0,4cm}

\section{The local Aligned Ordered Hilbert Scheme over a Grassmann Variety.} 

 \vspace{0,4cm}
 
 In the first part of this section, we assume that $\mathrm{Spec} \ R \simeq \mathbb{A} ^{2N-2}$ is an affine open set of the Grassmann variety  $G(1,N)$.
  \vspace{0,2cm}

We consider  an  affine line  $L \subset \mathbb{A}^N= \mathop{\rm Spec} {\CC}[x_1,...,x_{N-1},z] $,  with equations $x_1=...=x_{N-1}=0$.
  \vspace{0,2cm}

Let  $\mathop{\rm Spec} R=\mathop{\rm Spec} \ {\CC}[u_1,...,u_{N-1},v_1,...,v_{N-1}] $ be such that 

- the line $L$ corresponds to the origin $ (0,...,0) \in\mathop{\rm Spec} \ R$,

- the local system of parameters  $u_i, v_i$ (for $1\leq i \leq N-1$) of  $G(1,N)$ and the indeterminate $z$ parametrizing the canonical affine line over  $\mathop{\rm Spec}  R$ verify the relations 
$$x_i=u_iz+v_i, \mbox{ for } \  i=1,...,N-1.$$

 The canonical inclusions 
$${\CC}[x_1,...,x_{N-1},z]={\CC}[u_1z+v_1,...,u_{N-1}z+v_{N-1},z] \subset {\CC}[u_1,...,u_{N-1},v_1,...,v_{N-1},z] \supset {\CC}[u_1,...,u_{N-1},v_1,...,v_{N-1}].$$
induce the morphisms $\pi$  and $\psi$ in the  following commutative diagram: 
 \vspace{0,5cm}

$$ \xymatrix{
\mathop{\rm Spec} R \ar@{=}[d]   &\mathop{\rm Spec} R[z] \ar[l]  \ar@{=}[d]  \ar[r] & \mathbb{A}^N\ar@{=}[d] \\
 \mathop{\rm Spec}  {\CC}[u_1,...,u_{n-1},v_1,...,v_{n-1}]  \ar@{^{(}->}[d] & \mathop{\rm Spec}  {\CC}[u_1,...,u_{n-1},v_1,...,v_{n-1},z]  \ar[l]_\psi  \ar@{^{(}->}[d] \ar[r]^-\pi  & \mathop{\rm Spec}  {\CC}[x_1,...,x_{N-1},z]   \ar@{^{(}->}[d] \\  
G(1,N) &  {\mathcal  I} \ar[l]  \ar[r]&  {\PP}^N. \\
}$$

 \vspace{0,5cm}

In the affine space $\mathbb{A}^N=\mathop{\rm Spec} {\CC}[x_1,...,x_{N-1},z] $, we consider now an irreducible  complete intersection $Y$ 
and a finite  scheme $Z \subset Y\cap L$ supported in $r$ distinct ordered  points $(0,...,0,a_i) \in Y\cap L\subset \mathop{\rm Spec} {\CC}[x_1,...,x_{N-1},z]$ and with multiplicity  $k_i$  
at  $(0,...,0,a_i)$. We assume that  $Y$ is smooth at all points  of $Z$.

Let $c$ be the codimension of $Y$  in  $\mathbb{A}^N$.  We can find a system of $c$ hypersurfaces $G_i \subset \mathbb{A}^N$ 
 defined by polynomials $g_i \in {\CC}[x_1,...,x_{N-1},z]$  such that
 
$$ \cap G_i = Y,$$ 

$$Y\cap L = G_1 \cap L,$$ 

$$L \subset G_s \  \  2 \leq s \leq c.$$
\vspace{0,2cm}

We denote by $\varGamma = \pi^{-1}(Y) \subset \mathop{\rm Spec}  {\CC}[u_1,...,u_{n-1},v_1,...,v_{n-1},z]=\mathop{\rm Spec} R[z]$ the inverse image of $Y$. It is  cut out  by
the equations $g_s(u_iz+v_i,z)=0 \in R[z]$,  for $s=1,...,c$. To  the morphism
$$\phi : \varGamma = \mathop{\rm Spec}  R[z]/(g_1,...,g_c)  \rightarrow  \mathop{\rm Spec} R$$ 
is associated the ordered aligned Hilbert Scheme 
$$OH_{(k_1,...,k_r)}(\phi)= \cap_{s=1}^c OH_{(k_1,...,k_r)}(g_s) \subset \mathop{\rm Spec}  R[z_1,...,z_r].$$
We intend to study the tangent space of this Hilbert Scheme at the point  $\{Z\} \in OH_{(k_1,...,k_r)}(\phi)$. 
We recall that  the $r$ distinct points of $(0,...,0,a_i) \in Y\cap L\subset \mathop{\rm Spec} {\CC}[x_1,...,x_{N-1},z]$ are  smooth in $Y$. As a consequence, we note that  $\varGamma$ is smooth at  the $r$ points $(0,....,0,a_i) \in  \mathop{\rm Spec} {\CC}[u_i,v_j,z]$.

\vspace{0,3cm}

We begin with describing  the equations and the cotangential equations of $OH_{(k)}(g) \subset \mathrm{Spec} \ R[z_1]$ in a neighborhood of $(o,a)$ 
 for a polynomial  $g(z) \in {\CC}[x_i,z_1] $.  We write $z$ for $z_1$. The proof of the following lemma is straightforward (essentially contained in the statement).

\begin{lem}

Consider  $g(z) \in {\CC}[x_i,z]\ \mbox{ with } i=1,...,N-1$ and  $g=g(u_iz+v_i,z) \in R[z]$, with  $i=1,...,N-1$.
We assume that  
$ cl(g) \in (R/{\mathcal M})[z]$  has multiplicity  $\geq   k$  at  the  point   $a \in {\CC}$.

1) There exists a unique  decomposition

 $$g(z) \equiv  p(z)+\sum (u_iz+v_i)q_{i}(z) \quad  \mathop{\rm mod} \ {\mathcal M}^2R[z], \quad p,q_{i} \in {\CC}[z] , \quad    p\in (z-a)^{k}{\CC}[z],$$

 \vspace{0,3cm}
 
2) It induces   a decomposition (we recall the convention $q^{(-1)}(z)=0$) 

  $$g^{(s)}(z) \equiv  p^{(s)}(z) + \sum u_iq_i^{(s-1)}(z) +  \sum_i (u_iz+v_i)q_{i}^{(s)}(z)  \quad  \mathop{\rm mod}{\mathcal M}^2R[z]$$

and decompositions

 $$g^{(s)}(a)  \equiv \sum_i  u_iq_{i}^{(s-1)}(a) + \sum_i (u_ia+v_i)q_{i}^{(s)}(a)  \quad  \mathop{\rm mod}  {\mathcal M}^2R[z] ,  \quad s<k-1,$$

 $$g^{(k-1)}(a) + (z-a)g^{(k)}(a)   \equiv   $$
 $$(z-a)p^{(k)}(a)
+\sum_i  u_iq_{i}^{(k-2)}(a) + \sum_i (u_ia+v_i)q_{i}^{(k-1)}(a) 
 \quad   \mathop{\rm mod} \ ({\mathcal M},(z-a))^2.$$

\label{localG}

\end{lem}

Note here (once again) that if $g$ has multiplicity $>k$ at $a$, then $p^{(k)}(a) \in (z-a)$ and
$$ (z-a)p^{(k)}(a) \in ({\mathcal M},(z-a))^2.$$ 
This is why we  introduce the following unpleasant  convention (notation).
\vspace{0,3cm}

\emph{If $e_{j}$ is the multiplicity of $g_1$ at the point $a_j$, we put $h_j=k_j$ if $e_j>k_j$, and $h_j=k_j-1$ if $e_j=k_j$.}

\vspace{0,3cm}

In order to apply  Lemma \ref{localG} to the polynomials $g_s(z)$, we consider  the unique  decompositions

$$g_1(z) \equiv  p(z)+\sum (u_iz+v_i)q_{1,i}(z) \quad   \mathop{\rm mod}{\mathcal M}^2R[z], \quad p,q_{1,i} \in {\CC}[z], \quad \  p\in \cap_j (z-a_j)^{k_j}{\CC}[z],$$

$$ g_t(z) \equiv  \sum (u_iz+v_i)q_{t,i}(z)  \quad   \mathop{\rm mod} {\mathcal M}^2R[z], \quad q_{t,i} \in {\CC}[z], \quad t>1.$$

\vspace{0,2cm}

These decompositions  will play an important part in the proof of Theorem \ref{ohilbsch}. We choose to underline here the following Proposition which is a straightforward  consequence of Lemma \ref{hilbhyper} and Lemma \ref{localG}.

\begin{pro}
The aligned ordered Hilbert Scheme $OH_{(k_1,...,k_r)}(\phi)$ is  smooth   of expected dimension $2N-\nolinebreak 2+r-ck$ at the point 
$\{Z\}=(o,a_1,...,a_r)\in \mathop{\rm Spec} R[z_1,...,z_r]$, 
where $a_i\neq a_j$ for $i\neq j$,   if and only if 
 the following elements of ${\mathcal M}/{\mathcal M}^2$ are linearly independent.

$$ \sum_i  u_iq_{1,i}^{(s-1)}(a_j) + \sum_i (u_ia_j+v_i)q_{1,i}^{(s)}(a_j), \quad j=1,...,r,  \quad 0\leq  s\leq h_j-1,$$

$$ \sum_i  u_iq_{t,i}^{(s-1)}(a_j) + \sum_i (u_ia_j+v_i)q_{t,i}^{(s)}(a_j),  \quad t>1, \quad  j=1,...,r,  \quad 0\leq s 
\leq k_j-1.$$

\label{Grass}

\end{pro}

\vspace{0,5cm}

In the second part of this section $\mathop{\rm Spec} R \simeq \mathbb{A}^{N-1}$ is an affine open set of the Grassmann variety  ${\PP}^{N-1}(\beta)$ parametrizing  the lines of  ${\PP}^{N}$ through a point $\beta \in {\PP}^{N}$.
\vspace{0,2cm}

More precisely, from here we fix a point $\beta=(0,...,0,b)$, general in the line $ L \subset \mathbb{A}^N$. We recall that $ \mathop{\rm Spec} {\CC}[u_i,v_j]$ is an affine open set in $G(1,N)$.
The intersection of this open set with the closed subvariety  ${\PP}^{N-1}(\beta) \subset G(1,N)$ (the lines through $\beta$) is 
$$ \mathop{\rm Spec}  {\CC}[u_i,v_j]/(u_ib+v_i).$$

We put
 $$R_b =R/(u_i b+v_i)$$
and we denote by $\bar{u_i}$ and $\bar{v_i}$  the classes of $u_i$ and $v_i$ in $R_b$. The relations
 $\bar{v_i}=-b\bar{u_i}$ need no comment and $\bar{u_i}$ is a system of generators (regular parameters) of the maximal ideal
  $${\mathcal M}_b = {\mathcal M}/(u_i b+v_i) \subset R_b.$$

 In the inverse image of $ \mathop{\rm Spec}  R_b$, in the incidence variety point/line,  we consider, as earlier, the affine open set 
 $\mathop{\rm Spec} R_b[z]$.
 This is an open affine variety in the blowing-up   $ \tilde{{\PP}}^N $ of
$ {\PP}^{N}$  at the point $\beta$. We observe now  the following commutative diagram:

\vspace{0,3cm}

$$ \xymatrix{
\mathop{\rm Spec}  R_b  \ar@{=}[d]   &\mathop{\rm Spec}  R_b[z] \ar[l]  \ar@{=}[d]  \ar[r] & \mathbb{A}^N\ar@{=}[d] \\
 \mathop{\rm Spec} {\CC}[\bar{u_1},...,\bar{u}_{n-1}]  \ar@{^{(}->}[d] & \mathop{\rm Spec} {\CC}[\bar{u_1},...,\bar{u}_{n-1},z]  \ar[l]_{\psi_{\beta}}  \ar@{^{(}->}[d] \ar[r]^-{\pi_{\beta}}  & \mathop{\rm Spec}  {\CC}[x_1,...,x_{N-1},z]   \ar@{^{(}->}[d] \\  
{\PP}^{N-1} (\beta) &\tilde{{\PP}}^{N}   \ar[l]  \ar[r] &  {\PP}^N. \\
}$$

\vspace{0,4cm}

We note that $\pi_{\beta}$  is the blowing-up of the point $(0,...,0,b)\in  \mathop{\rm Spec} {\CC}[x_1,...,x_{N-1},z] $.
 
 We recall that $L\subset \mathbb{A}^N$ is the affine line with equations $x_i=0  \ \mbox{ with } i=1,...,N-1$. Its inverse image in 
 $\mathop{\rm Spec} R_b[z] $ is  cut out by the equations 
$$0=x_i= \bar{u_i}z + \bar{v_i} = \bar{u_i}(z-b), \quad  i=1,...,N-1.$$

We also recall  that $Y \subset \mathbb{A}^N$ is the complete intersection
of $c$ hypersurfaces $G_i \subset \mathbb{A}^N$ 
 defined by polynomials $g_i \in {\CC}[x_1,...,x_{N-1},z]$  such that

$$Y\cap L = G_1 \cap L, $$ 

$$L \subset G_s, \quad \  2 \leq s \leq c.$$
\vspace{0,2cm}

The inverse image (proper transform) 
$\varGamma_b= \pi_{\beta}^{-1}(Y) \subset \mathop{\rm Spec} R_b[z] $ 
is cut out by the $c$ equations 
$$g_s(\bar{u_i}(z-b),z)=0, \quad s=1,...,c.$$

We put $g_{s,b} =  g_s(\bar{u_i}(z-b),z) \in R_b[z]$, denote $R_b [z]/(g_{1,b},...,g_{c,b}) =R[z]/(g_1,...,g_c)  \otimes_R R_b$ and  
consider the morphism 
$$\phi_b :  \varGamma_b= \mathrm{Spec} \ R_b [z]/(g_{1,b},...,g_{c,b}) \subset  \mathop{\rm Spec}R_b[z] \rightarrow  \mathop{\rm Spec} R_b.$$

We intend to study the tangent space of the Hilbert Scheme at the point  
$$\{Z\} \in OH_{(k_1,...,k_r)}(\phi_b)= \cap_{s=1}^c OH_{(k_1,...,k_r)}(g_{s,b}) \subset \mathrm{Spec} \ R_b[z_1,...,z_r].$$ 
Let us recall that  the $r$ distinct points of $(0,...,0,a_i) \in Y\cap L\subset \mathop{\rm Spec} {\CC}[x_1,...,x_{N-1},z]$ are  smooth in $Y$. As a consequence, we note that  $\varGamma_b$ is smooth at  the $r$ points $(0,....,0,a_i) \in  \mathop{\rm Spec} {\CC}[\bar{u_i},z]$.

\vspace{0,3cm}

The cartesian diagram

$$\begin{CD}
\varGamma_b @>\phi_b >>  \mathop{\rm Spec}R_b\\
@VVV @VVV\\
\varGamma  @>\phi>>  \mathop{\rm Spec} R
\end{CD}$$
\vspace{0,2cm}
and Proposition \ref{defhilb}  3) imply the following

\begin{pro}
$$OH_{(k_1,...,k_r)}(\phi_b)=OH_{(k_1,...,k_r)}(\phi)\times_ {\mathop{\rm Spec} R} ( \mathop{\rm Spec} R_b).$$

\label{hilbuniv}

\end{pro}

\vspace{0,3cm}

Next  we intend to give explicit necessary and sufficient conditions  for the smoothness of $OH_{(k_1,...,k_r)}(\phi_b)$ at $(o,a_1,...,a_r)$.
As in the preceding case, we begin by describing the local equations (in a neighborhood of a point $(0,a)$) of  the Hilbert Scheme $OH_{(k)}(g_b) \subset   \mathrm{Spec} \ R_b[z]$, where $g_b$ is of the forme $g_b(z)=g(\bar{u_i}(z-b),z)$, with 
$g\in ( {\CC}[x_1,...,x_{N-1},z]$.  Lemma \ref{localG}   specializes immediately in the following way:

\begin{lem} 

\vspace{0,3cm}

Consider a polynomial $g(x_i,z) \in {\CC}[x_i,z]$,  with  $ i=1,...,N-1$, and  
$ g_b=g(\bar{u_i}(z-b),z) \in R_b[z], \mbox{ with } i=1,...,N-1.$
  Assume that the polynomial 
$ cl(g_b) \in (R_b/{\mathcal M}_b)[z]$  has multiplicity  $\geq   k$  at  the  point    $a \in {\CC}$.

The unique decomposition (as before,  we follow the convention  $q_i^{(-1)}=0$)

 $$g_b(z) \equiv  p(z)+\sum \bar{u_i}(z-b)q_{i}(z)\quad  \mathop{\rm mod}  {\mathcal M}^2R[z], \quad  p,q_i \in {\CC}[z],  \quad  p\in (z-a)^{k}{\CC}[z],$$
induces, for all $s$,  a decomposition  
 
 $$g_b^{(s)}(z)  \equiv p^{(s)}(z) + \sum_i  \bar{u_i}q_{i}^{(s-1)}(z) + \sum_i \bar{u_i}(z-b)q_{i}^{(s)}(z)  \quad    \mathop{\rm mod}  {\mathcal M}_b^2R_b[z].$$
and  decompositions

 $$g_b^{(s)}(a)  \equiv  \sum_i  \bar{u_i}q_{i}^{(s-1)}(a) + \sum_i \bar{u_i}(a-b)q_{i}^{(s)}(a) = $$
 $$\sum_i  \bar{u_i}[q_{i}^{(s-1)}(a) + (a-b)q_{i}^{(s)}(a)]
 \quad  \mathop{\rm mod}  {\mathcal M}_b^2R_b[z],
 \quad   s<k-1,$$

$$g_b^{(k-1)}(a) + (z-a)g_b^{(k)}(a)   \equiv $$
$$(z-a)p^{(k)}(a)
+\sum_i  \bar{u_i}q_{i}^{(k-2)}(a) + \sum_i \bar{u_i}(a-b)q_{i}^{(k-1)}(a)
 \quad    \mathop{\rm mod} ({\mathcal M}_b,(z-a))^2.$$
\end{lem} 

We recall that that $ g_{t,b}(z) \in {\mathcal M}_bR_b[z]$ for $ t>1$.  As an immediate consequence, we get the following result (to be compared with Proposition \ref{Grass}):

\begin{pro}
The aligned ordered Hilbert Scheme $OH_{(k_1,...,k_r)}(\phi_b)$ is  smooth of expected dimension $N-1+r-kc$ at the point 
$\{Z\}=(o,a_1,...,a_r) \in \mathop{\rm Spec} R[z_1,...,z_r]$  with $a_i\neq a_j$ for $i\neq j$  if and only if   the following elements of ${\mathcal M}_b/{\mathcal M}_b^2$ are  linearly independent:

$$ \sum_i  \bar{u_i}[q_{1,i}^{(s-1)}(a_j) + (a_j-b)q_{1,i}^{(s)}(a_j)], \quad j=1,...,r, \quad  0\leq  s\leq h_j-1,$$

$$ \sum_i  \bar{u_i}[q_{t,i}^{(s-1)}(a_j) + (a_j-b)q_{t,i}^{(s)}(a_j)], \quad  t>1, \quad  j=1,...,r , \quad 0\leq s 
\leq k_j-1.$$

\label{grassloc}

\end{pro}

\qed

Finally in this section,  we observe that, for $t$ and $j$ fixed and $b\neq a_j$, the vector subspaces
of ${\mathcal M}_b/{\mathcal M}_b^2$, generated by 

$$ \sum_i  \bar{u_i}[q_{1,i}^{(s-1)}(a_j) + (a_j-b)q_{1,i}^{(s)}(a_j)]), \quad  0\leq s 
\leq h_j-1,   \quad h_j-1,$$
and 
$$ \sum_i  \bar{u_i}[q_{t,i}^{(s-1)}(a_j) + (a_j-b)q_{t,i}^{(s)}(a_j)]), \quad  0\leq s 
\leq k_j-1,  \quad t>1, $$
on the one hand, and 
$$\sum_i  \bar{u_i}q_{1,i}^{(s)}(a_j), \quad  0\leq s 
\leq h_j-1, $$
and 
$$\sum_i  \bar{u_i}q_{t,i}^{(s)}(a_j), \quad  0\leq s 
\leq k_j-1,  \quad t>1, $$
on the other hand,  coincide with each other.  This proves  the equivalence $1) \Leftrightarrow  2) $ in the following Corollary (of the previous Proposition):

\begin{coro}
If  $b\neq a_j$ for all $j$, the following equivalent conditions are equivalent.

1) The aligned ordered Hilbert Scheme $OH_{(k_1,...,k_r)}(\phi_b)$ is smooth of expected dimension $N-1+r-ck$  at the point 
$\{Z\}=(o,a_1,...,a_r)$.

2) The following elements of ${\mathcal M}_b/{\mathcal M}_b^2$ are  linearly independent:

$$ \sum_i  \bar{u_i}q_{1,i}^{(s)}(a_j), \quad j=1,...,r, \quad  0\leq  s\leq h_j-1,$$

$$\sum_i  \bar{u_i}q_{t,i}^{(s)}(a_j),  \quad  t>1, \quad  j=1,...,r, \quad 0\leq s 
\leq k_j-1.$$

3) The  matrix $(q_{t,i}^{(s)}(a_j))$  with $N-1$ rows and $(c-1)k+\sum_j h_j$ columns has   maximal rank.

\label{linearalgebra}

\end{coro}

\proof  We have seen  $1) \Leftrightarrow  2) $.  The equivalence $2) \Leftrightarrow  3) $ is an obvious consequence of the fact that   
 $(\bar{u_i})_i$ is a regular system of generators of ${\mathcal M}_b$. 
\qed

\vspace{0,2cm}

Note to conclude this section that condition 3) does not depend on $b$ (this will be a crucial point in the proof (by induction)  of Theorem \ref{ohilbsch}). In other words, if $\{Z\}=(o,a_1,...,a_r)$ is a smooth point of 
$OH_{(k_1,...,k_r)}(\phi_b)$ and if $\beta' =(0,...,0,b') \in L$, with $b'\neq a_i, \quad i=1,...,r$, then $\{Z\}=(o,a_1,...,a_r)$ is also a smooth point of 
$OH_{(k_1,...,k_r)}(\phi_{b'})$.

\vspace{0,3cm}

\section{Proof of Theorem \ref{ohilbsch}}

\vspace{0,6cm}

In the previous section, we have studied the configuration of a line $L$, a quasi-projective complete intersection $Y \subset {\PP}^N$ and a general point $\beta \in L$. We studied a  finite scheme 
$Z \subset L\cap Y$, with support in the smooth locus of $Y$ and with multiplicities $(k_1,...,k_r)$ in $r$ distinct points of $L\cap Y$. This is a point of $OH_{(k_1,...,k_r)}(Y)$.   We recall that the inverse images $\varGamma$ and $\varGamma_b$ of $Y$, in the incidence varieties $ {\mathcal  I}$ and ${\PP}^N(\beta)$, fit in
 a cartesian diagram

$$\begin{CD}
\varGamma_b @>\phi_b >>  \mathrm{Spec} \ R_b\\
@VVV @VVV\\
\varGamma @>\phi>>  \mathrm{Spec} \ R
\end{CD}$$\vspace{0,2cm}

We keep these notations in mind and we come back to the composed projection morphism described in Theorem \ref{ohilbsch}

$$\theta_{(k_1,...,k_r)} : OH_{(k_1,...,k_r)}(Y) \times_{G}  {\mathcal  I} \rightarrow {\mathcal  I}  \rightarrow  {\PP}^{N}.$$

In order to study this morphism in a neighorhood of  the locally closed subscheme
$$\{Z\}\times_{\mathop{\rm Spec} R} \mathop{\rm Spec} R[z] =(o,a_1,...,a_r) \times_{\mathop{\rm Spec} R} \mathop{\rm Spec} R[z] \subset  \mathop{\rm Spec} R[z_1,...,z_r,z] ,$$ we observe the following  commutative diagram (where all up vertical arrows are closed immersions):

\vspace{0,3cm}

$$ \xymatrix{
 \mathop{\rm Spec} R[z_1,...,z_r]    &\mathop{\rm Spec}R[z_1,...,z_r][z]  \ar[l]    \ar[r] &   \mathop{\rm Spec}R[z]  \ar@{=}[d]  \ar[r]    &\mathbb{A}^N\ar@{=}[d] \\
OH_{(k_1,...,k_r)}(\phi)  \ar[u]  \ar@{^{(}->}[d] &OH_{(k_1,...,k_r)}(\phi) \times_{\mathop{\rm Spec} R} (\mathop{\rm Spec}R[z] ) \ar[u] \ar[l]  \ar@{^{(}->}[d] \ar[r] &
\mathop{\rm Spec} R[z]   \ar@{^{(}->}[d] \ar[r]^-p &\mathbb{A}^N \ar@{^{(}->}[d]\\  
 OH_{(k_1,...,k_r)}(Y) &OH_{(k_1,...,k_r)}(Y) \times_{G}  {\mathcal  I}   \ar@{=}[d] \ar[l] \ar[r] & {\mathcal  I}  \ar[r]    &  {\PP}^N  \ar@{=}[d]\\
&      OH_{(k_1,...,k_r)}(Y) \times_{G}  {\mathcal  I}  \ar[r] &\theta_{(k_1,...,k_r)}\ar[r] &{\PP}^N. \\
}$$

\vspace{0,4cm}

We recall here that the equations of $L$ in $\mathbb{A}^N = \mathop{\rm Spec}{\CC}[x_1,...,x_{N-1},z]$
 are $x_1=...=x_{N-1}=0$ and that $\beta =(0,...,0,b) \in L$.

\begin{pro} 

$\theta_{(k_1,...,k_r)}^{-1}(\beta) \cap (OH_{(k_1,...,k_r)}(\phi) \times_{\mathop{\rm Spec} R} (\mathop{\rm Spec} R[z])) =OH_{(k_1,...,k_r)}(\phi_b)$.
\label{fibre}
\end{pro}

\proof

We begin with describing  the fiber  $p^{-1}(\beta) \subset \mathop{\rm Spec}R[z]$. The maximal ideal of ${\CC}[x_1,...,x_{N-1},z ]$ corresponding to $\beta$ is $(x_1,...,x_{N-1},z-b)$ and
$$R[z]/(x_1,...,x_{N-1},z-b)=R[z]/(u_iz+v_i,z-b)=R[z]/(u_ib+v_i,z-b)=R_b.$$

But we have seen (proposition \ref{hilbuniv}) that 
$$OH_{(k_1,...,k_r)}(\phi_b)=OH_{(k_1,...,k_r)}(\phi)\times_ {\mathop{\rm Spec} R} ( \mathop{\rm Spec} R_b),$$
so the following commutative diagram  proves   our Proposition:

\vspace{0,3cm}

$$\begin{CD}
         \mathop{\rm Spec} R[z_1,...,z_r]@<<<  \mathop{\rm Spec}R[z_1,...,z_r][z] @>>>  \mathop{\rm Spec} R[z]@>>> \mathbb{A}^N\\
                              @AAA                                       @AAA                                                              @|                                   @|\\    
         OH_{(k_1,...,k_r)}(\phi) @<<<  OH_{(k_1,...,k_r)}(\phi) \times_{\mathop{\rm Spec}R} (\mathop{\rm Spec}R[z])  @>>>  \mathop{\rm Spec}R[z] @>>> \mathbb{A}_N \\
                              @|                                            @AAA                        @AAA                                                    @AAA          \\                                 
      OH_{(k_1,...,k_r)}(\phi) @<<<  OH_{(k_1,...,k_r)}(\phi) \times_{\mathop{\rm Spec}R} (\mathop{\rm Spec} (R[z]\otimes {\CC}(\beta)))   @>>>  \mathop{\rm Spec} (R[z]\otimes {\CC}(\beta)) @>>>  \beta \\                    
                    @| @| @| @| \\                       
            OH_{(k_1,...,k_r)}(\phi) @<<<  OH_{(k_1,...,k_r)}(\phi)\times_{\mathop{\rm Spec}R} \mathop{\rm Spec} R_b @>>>  \mathop{\rm Spec}R_b @>>>  \beta\\
         @AAA                                           @|                               @.            @.\\
         OH_{(k_1,...,k_r)}(\phi_b) @= OH_{(k_1,...,k_r)}(\phi_b). @.            @.
      \end{CD}$$

\qed

\vspace{0,4cm} 

We can now proceed with the proof, by induction on $k$, of  Theorem \ref{ohilbsch} (which we recall).

\begin{theo} (Aligned Ordered Hilbert Scheme Theorem)   Let  $X$ be a smooth connected dimension $n$ quasi-projective variety embedded in  ${\PP}^{N}$, with $N=n+c$.

For $k=k_1+...+k_r$, with $k_i>0$, let $OH_{(k_1,...,k_r)}(X) $ be the  ordered Hilbert scheme  of aligned,  finite, degree $k$ subschemes of $X$,
with (ordered) multiplicities $k_i$ at the ordered  points $x_i$ (possibly coinciding).  Consider  the natural projective line bundle 
$H_{(k_1,...,k_r)}(X) \times_{G}  {\mathcal  I} $   over
 $H_{(k_1,...,k_r)}(X)$ and the projection 
 $$\theta_{(k_1,...,k_r)} : OH_{(k_1,...,k_r)}(X) \times_{G}  {\mathcal  I} \rightarrow {\mathcal  I}  \rightarrow  {\PP}^{N}.$$
 
 The general fiber of  $\theta_{(k_1,...,k_r)}  $ is smooth of dimension $N-1+r-kc$.

\end{theo}

We apply the results of the preceding section in the the case $Y=X$ and we claim that Theorem \ref{ohilbsch} is a consequence of the next  proposition. We go on considering a point   $\{Z\}=(L,a_1,...,a_r)\in OH_{(k_1,...,k_r)}(X) $ corresponding to a finite scheme $Z\subset L\cap X$ with multiplicities $(k_1,...,k_r)$ at the distinct points 
$$(a_1,...,a_r) \in \mathop{\rm Spec} {\CC}[z] \cap X \subset  L\cap X.$$

\vspace{0,2cm}

\begin{pro}  Consider  $\{Z\}=(L,a_1,...,a_r) \in OH_{(k_1,...,k_r)}(X) $, 
with $k=\sum k_i >1$.

For $k_r=1$, define $\{Z'\}= (L,a_1,...,a_{r-1}) \in OH_{(k_1,...,k_{r-1},k_{r-1})}(X) $, where $Z' \subset Z$ is the finite, degree $k-1$,  subscheme of $Z$ with multiplicity $k_i$ at $a_i$ for $i\leq r-1$ and multiplicity $0=k_r-1$ at $a_r$.

For $k_r>1$, define $\{Z'\}= (L,a_1,...,a_{r}) \in OH_{(k_1,...,k_r-1)}(X )$, where $Z' \subset Z$ is the finite, degree $k-1$, subscheme of $Z$ with multiplicity $k_i$ at $a_i$ for $i\leq r-1$ and multiplicity $k_r-1$  at $a_r $.

\vspace{0,2cm}

If  $OH_{(k_1,...,k_r)}(X) $ is not smooth of dimension $2N-2+r-kc$  at $\{Z\}$ , then 

- for $k_r=1$, the point  $(\{Z'\},\beta) \in OH_{(k_1,...,k_{r-1})}(X) \times_{G}  {\mathcal  I} $ is a point of ramification for $\theta_{(k_1,...,k_{r-1})} $,  

- for $k_r>1$,   the point  $(\{Z'\},\beta) \in  OH_{(k_1,...,k_r-1)}(X)\times_{G}  {\mathcal  I}  $  is a point of ramification for $\theta_{(k_1,...,k_r-1)}. $

\label{indprincilocal}

\end{pro}

 \textit{Proof}  of Theorem  \ref{ohilbsch}. 

We assume that the proposition is true and we proceed by induction on $k$.
Note that for $k =1$  the aligned ordered Hilbert scheme $OH_{(1)}(f) $ is smooth, hence 
$OH_{(1)}(X) \times_{G}  {\mathcal  I} $ is smooth and the general fiber of  $\phi_{(1)}  $ is smooth of dimension $N-1+1-c=n$ by  Bertini's Theorem.

Let $k$ be minimum number for which there exists a  partition $k=k_1+...+k_r, \quad k_i>0$  and such that the generic  fiber of  $\phi_{(k_1,...,k_r)}  $ fails to be smooth of dimension $N-1+r-kc$.  By Bertini's theorem, this implies that the inverse image of  the singular locus of $OH_{(k_1,...,k_{r})}(X)$ in $OH_{(k_1,...,k_{r})}(X) \times_{G}  {\mathcal  I} $ dominates ${\PP}^{N}$. Applying  proposition \ref{indprincilocal}, we find that the ramification locus of $\theta_{(k_1,...,k_{r-1})} $ (or $\theta_{(k_1,...,k_r-1)} $ if $k_r>1$) dominates ${\PP}^{N}$. This contradicts the minimality of $k$. 

\qed

\textit{Proof}  of Proposition \ref{indprincilocal}. 

Note that by  Proposition~\ref{lpointsdistinct} we can assume $a_i \neq a_j$ for $i\neq j$ (this is a key point!).
\vspace{0,3cm}

If  $OH_{(k_1,...,k_r)}(\phi) $ is not smooth of dimension $2N-2+r-kc$  at  $x=(L,a_1,...,a_r)$, then, by Proposition \ref{Grass},
 the elements

$$ \sum_i  u_iq_{1,i}^{(s-1)}(a_j) + \sum_i (u_ia_j+v_i)q_{1,i}^{(s)}(a_j), \quad  j=1,...,r, \quad 0\leq  s\leq h_j-1,$$
and 

$$\sum_i  u_iq_{t,i}^{(s-1)}(a_j) + \sum_i (u_ia_j+v_i)q_{t,i}^{(s)}(a_j), \quad  t>1, \ \  j=1,...,r , \quad  0\leq s 
\leq k_j-1$$
are linearly dependent in ${\mathcal M}/{\mathcal M}^2$.

Specializing in   ${\mathcal M}_b={\mathcal M}/(u_ib+v_i)$, we see that  the elements 

$$ \sum_i  \bar{u_i}[q_{1,i}^{(s-1)}(a_j) + (a_j-b)q_{1,i}^{(s)}(a_j)], \quad  j=1,...,r, \quad   0\leq  s\leq h_j-1,$$

$$\sum_i  \bar{u_i}[q_{t,i}^{(s-1)}(a_j) + (a_j-b)q_{t,i}^{(s)}(a_j)], \quad  t>1, \quad  j=1,...,r, \quad  0\leq s 
\leq k_j-1$$
are linearly dependent in ${\mathcal M}_b/{\mathcal M}_b^2$.
\vspace{0,2cm}

 In the special case $b=a_r$, we find that the elements    

$$1) \quad \sum_i  \bar{u_i}[q_{1,i}^{(s-1)}(a_j) + (a_j-a_r)q_{1,i}^{(s)}(a_j)], \quad   j=1,...,r-1, \quad  0\leq  s\leq h_j-1,$$
or equivalently 
$$\quad \sum_i  \bar{u_i}[q_{1,i}^{(s)}(a_j)], \quad   j=1,...,r-1, \quad  0\leq  s\leq h_j-1,$$

 $$2) \quad \sum_i  \bar{u_i}[q_{1,i}^{(s-1)}(a_r)], \quad 0\leq  s\leq h_r-1,$$
or equivalently
$$ \sum_i  \bar{u_i}[q_{1,i}^{(s)}(a_r)], \quad   \  0\leq  s\leq h_r-2=(h_r-1)-1,$$

$$3) \sum_i  \bar{u_i}[q_{t,i}^{(s-1)}(a_j) + (a_j-a_r)q_{t,i}^{(s)}(a_j)], \quad  t>1, \quad   j=1,...,r-1, \quad  0\leq s 
\leq k_j-1,$$
or equivalently
$$ \sum_i  \bar{u_i}[q_{t,i}^{(s)}(a_j)], \quad  t>1, \quad   j=1,...,r-1, \quad  0\leq s 
\leq k_j-1,$$

and 
$$4) \quad  \sum_i  \bar{u_i}q_{t,i}^{(s-1)}(a_r), \quad  t>1, \quad  0\leq s \leq k_j-1,$$
or equivalentely
$$ \sum_i  \bar{u_i}q_{t,i}^{(s)}(a_r), \quad   t>1, \quad  0\leq s  \leq k_r-2=(k_r-1)-1$$
are linearly dependent in ${\mathcal M}_{a_r}/{\mathcal M}_{a_r}^2$. 
\vspace{0,2cm}

Using then   Corollary~\ref{linearalgebra} and Proposition \ref{fibre}, one sees easily that 
\vspace{0,2cm}

 $$OH_{(k_1,...,k_{r-1})}(\phi_b)=\theta_{(k_1,...,k_{r-1})}^{-1}(\beta) \cap (OH_{(k_1,...,k_{r-1})}(\phi) \times_{\mathrm{Spec} \ R} (\mathrm{Spec} \ R[z]))$$
is singular at $\{Z'\}$ when $k_r=1$, and that

$$OH_{(k_1,...,k_{r}-1)}(\phi_b)=\theta_{(k_1,...,k_r-1)}^{-1}(\beta) \cap (OH_{(k_1,...,k_r-1)}(\phi) \times_{\mathrm{Spec} \ R} (\mathrm{Spec} \ R[z]))$$ is singular at $\{Z'\}$ when $k_r>1$.

\qed

\vspace{0,2cm}

As we already remarked, Proposition \ref{indprincilocal} implies  Theorem \ref{ohilbsch} which in turn implies 
 Theorem   \ref{hilbsch} which yields Theorem \ref{genproj}.

\vspace{0,3cm}

\section{Examples, questions and conjectures} 

 \textsc{Examples}
 
 \vspace{0,4cm}
 
\begin{exa}  As a first example, consider  a projected \textsc{Veronese surface} $X \subset {\PP}^4$ (yes projected in $ {\PP}^4$), and  a general projection  
 $X \rightarrow X_1 \subset  {\PP}^3$.
 \end{exa}
 
The Steiner surface $X_1$ is well known.  We describe its singular locus.

- $X_2=X_{\{1,1\}}$ is composed of  three lines  through a point  $x$  and  not in a plane.  

The normalization $\tilde{X}_2$ of $X_2$,  a fiber of the map  
$\phi_{\{1,1\}} : H_{\{1,1\}}(X) \times_{G}  {\mathcal  I} \rightarrow {\mathcal  I}  \rightarrow  {\PP}^{4},$
 is a disjoint union of three lines.

- $X_2$ has a closed subscheme $X_{\{2\}}$ composed of $6$ distinct pinch points, $2$ on each of the $3$ lines.

- The degree $1$ finite scheme  $X_3=X_{\{1,1,1\}} = \{x\}$ is the  triple locus  of $X_1$, as well as the singular and triple locus of  $X_2$.   We note that, as stated in Theorem \ref{genproj},  we have $X_{\{2,1\}}=\emptyset$, in other words $X_{\{2\}}$  and $X_3$ are disjoint.

\vspace{0,2cm}

\vspace{0,2cm}

\begin{exa} The Veronese surface $X \subset {\PP}^4$ is  one of the four  \textsc{Severi Varieties}. According to a celebrated result of F. Zak (see \ref{Zak}), if $X^n\subset {\PP}^N$ is a nondegenerate, 
dimension $n $, 
 smooth variety  with $N \leq 3n/2+1$, 
then $X$ is linearly complete except for   the four  projected Severi varieties, for which $n=2^k$  with $k=1,2,3,4$ and $N =3n/2+1$. 

\label{sevar}
 \end{exa}

We consider a projected  Severi variety  $X \subset {\PP}^{(3n/2)+1}$  and we  describe the singularities of a general projection  $X \rightarrow X_1 \subset {\PP}^{3n/2}.$ 

- $X_2=X_{\{1,1\}}$ is composed of  three ${\PP}^{n/2}$ through a point  $x$  and  not in a hyperplane. 
Its normalization $\tilde{X}_{\{1,1\}}$ is a disjoint union of three ${\PP}^{n/2}$.  We recall that  $\tilde{X}_{\{1,1\}}$ is a general  fiber of the map  $\phi_{\{1,1\}} : H_{\{1,1\}}(X) \times_{G}  {\mathcal  I} \rightarrow {\mathcal  I}  \rightarrow  {\PP}^{N}.$

- $X_2$ has a closed subscheme $X_{\{2\}}$ composed of three disjoint quadrics of dimension  $(n/2)-1$, one in each  ${\PP}^{n/2}$.

- The degree $1$ finite  scheme  $X_3=X_{\{1,1,1\}} = \{x\}$ is the  triple locus  of $X_1$  as well as the singular and triple locus of  $X_2$.   

- We note once again that the degree $6$ finite scheme $X_{\{2\}}$  is smooth and disjoint from $X_3$. Indeed $X_{\{2,1\}}=\emptyset$, as stated in Theorem \ref{genproj}.

\vspace{0,2cm}

   \vspace{0,2cm}
   
\begin{exa} Consider a general skew-symmetric map $6O_{{\PP}^3} (-1)\rightarrow  6O^{{\PP}_3}$. The cubic surface defined by its degree $3$ pfaffian is smooth and equipped with a projective ${\PP}^1$-bundle smoothly embedded in ${\PP}^5$ as a $3$-fold of degree $7$, well known as the  \textsc{Palatini} $3$-fold.

\label{palvar}

\end{exa}

 We describe now the singularities of a general projection $X \rightarrow X_1 \subset {\PP}^4$ of a Palatini $3$-fold.
 
- $X_2=X_{\{1,1\}}$ is an irreducible surface of degree $11$ whose smooth normalization $\tilde{X}_{\{1,1\}}$ is a fiber of $\phi_{\{1,1\}} : H_{\{1,1\}}(X) \times_{G}  {\mathcal  I} \rightarrow {\mathcal  I}  \rightarrow  {\PP}^{5}.$

- $X_3=X_{\{1,1,1\}} $ is the singular and the triple locus of $X_2$. It is composed of four lines through a point $x$ and generating  ${\PP}^4$. 
The normalization  of  $X_{\{1,1,1\}} $ is a disjoint union of four lines.

- $X_2$ contains a pinch curve $X_{\{2\}}$ of degree $22$. 

- The scheme  $X_{\{2\}} \cap X_3=X_{\{2,1\}}$ is the singular locus of $X_{\{2\}}$. It is composed of $24$ distinct points, six on each of the four lines.

- $X_4= \{x\}$ is a degree $1$ finite scheme. By Theorem \ref{genproj} we have  $X_4\cap X_{\{2,1\}}=X_{\{2,1,1\}}=\emptyset$.

   \vspace{0,4cm}
   
\begin{exa} Consider an elliptic quintic ruled surface $S \subset {\PP}^4$ (the lines of $S$ are parametrized by a section of $G(1,4)$ by a general ${\PP}^4$ in the Plücker space).
\end{exa}

We describe the singularities of a general projection $S \rightarrow S_1 \subset {\PP}^3$.

   - The double locus $S_{\{1,1\}}=S_2 \subset S_1$ is a smooth quintic elliptic curve.
   
   - The triple locus $S_3=S_{\{1,1,1\}} $ is empty. This deserves a comment, see  Theorem \ref{Aure}.
   
   - There are ten  distinct pinch points on $S_2$, i.e. $S_{\{2\}}$ is a smooth, degree $10$, finite scheme.
  
    \vspace{0,4cm}
    
 \textsc{Questions and conjectures} 
    \vspace{0,4cm}
    
Our first and main question is classical.

\textit{Let $X \subset {\PP}^N$ be a dimension $n$  smooth variety, not contained in a hypersurface of degre $<k$.
For which $(n,N,k)$ do the $k$-secant lines to   $X $ fill up the space?}
       \vspace{0,2cm}
       
       We know by Theorem \ref{hilbsch} that $k(N-n-1)\leq N-1$ is a necessary consequence.
          \vspace{0,2cm}

For $k=2$ the complete answer was given by F. Zak (see for example \cite{9} or \cite{6}).

 \begin{theo} (F. Zak) 1) If $N-1-3(N-n-1) \geq -2$, then the $2$-secant lines to a nondegenerate dimension $n$ smooth variety $X\subset {\PP}^N$ fill up the ambient space except for  the four (nonprojected)  Severi Varieties, in which case  $n=2^k$ with $k=1,2,3,4$ and $N-1-3(N-n-1) = -2$ (i.e. $N =3n/2+2=3.2^{k-1}+2$).  

 \label{Zak}
 
  \end{theo}
  
  It is not irrelevant to recall that the $2$-secant lines to a Severi variety $X$ fill up a cubic hypersurface of ${\PP}^N$.
        \vspace{0,1cm}
  
  Note also  that the nonprojected Severi Varieties are   cut out by quadric hypersurfaces.
  
      \vspace{0,3cm}
  
  For $k=3$, the question is open except for $N=4$, in which case A. Aure  proved the following result (see \cite {2}):
  
   \begin{theo} (A. Aure)  Elliptic quintic scrolls in ${\PP}^4$ are the only smooth surfaces not contained in a quadric hypersurface  whose $3$-secant lines do not fill up ${\PP}^4$.
   
   \label{Aure}

  \end{theo}
  
  The $3$-secant lines to a quintic elliptic surface fill up a quintic hypersurface of  ${\PP}^4$.
          \vspace{0,1cm}
          
   Note that elliptic quintic scrolls in ${\PP}^4$ are cut out by cubic hypersurfaces.

   \vspace{0,2cm}
   
This suggests

    \begin{conj}  There exists  a function $\phi (k)$ such that  for any dimension $n$,smooth variety  $X \subset  {\PP}^N$ not contained in a hypersurface of degree $<k$ sone has: 
    
 1) if  $N-1-(k+1)(N-n-1) > \phi (k)$, then  the $k$-secant lines to $X$ fill up the ambient space;
   
 2) if $N-1-(k+1)(N-n-1) = \phi (k)$ and the $k$-secant lines to $X$ do not fill up the ambient space, then  $X$ is cut out by hypersurfaces of degree $k$.
        
   \end{conj}
   
   From Zak's Theorem we  get $\phi (2)=-2$ and basing on  Aure's Theorem we conjecture $\phi (3)=-1$.

   \vspace{0,2cm}

Our second series of questions concerns the irreducibility of the loci $X_k$ of a general projection of a smooth variety.

We begin with recalling  Franchetta's famous theorem.

     \begin{theo} (A. Franchetta)  The Veronese surface in ${\PP}^4$ is the only smooth projective surface  whose general projection to ${\PP}_3$ has a  reducible double locus.
  \end{theo}
  
  Of course, Franchetta does not assume  $S \subset {\PP}^4$,   but  by Bertini's Theorem  this is the only difficult case.
  
  As a comment to this result, we recall that 
  
  - if $X$ is a projected  Severi variety $X^n \subset {\PP}^{(3n/2)+1}$, then the locus $X_2$ is a 
  union of three ${\PP}^{n/2}$ through a point (Example \ref{sevar}),
  
  - if $X$ is a Palatini $3$-fold, then the locus $X_3$ is a  union of four  lines through a point (Example \ref{palvar}),
  
  - there exists a dimension $6$ smooth variety $X\subset {\PP}^9$  whose general projection has a reducible triple locus. $X_3$  is a union of four planes passing through a  point $x$, with  $X_4= \{x\}$).
  
  \vspace{0,3cm}

 We dare a bold conjecture
 
   \begin{conj}  Let $X \subset {\PP}^N$ be a dimension $n$, smooth irreducible variety.
   
  1)  If $(k +1)(N-n-1) < (N-1)$ then the locus $X_k$ of a general projection of $X$ is irreducible.
  
  2) If $(k +1)(N-n-1) = (N-1)$ and  the locus $X_k$ of a general projection of $X$ is reducible, then the finite scheme $X_{k+1}$ has degree $1$ and
  $X_k$ is a union of $k+1$ linear space ${\PP}^{N-n-1}$ passing through the point of $X_{k+1}$.

 \end{conj}

Our conjecture is  related to the following  conjecture of F. Zak (\cite{10}):

 \begin{conj} Let $X^n \subset  {\PP}^N$ with  $N-1\geq (k+1)(N-n-1)$ be a  nondegenerate irreducible variety (non necessarily smooth). Consider a general projection $X \rightarrow {\PP}^{N-1}$. Then the locus  $X_k \subset {\PP}^{N-1}$  of points whose fiber has degree $\geq k$) is connected. 
   \end{conj}

   To conclude this paper, we note here a theorem, a related conjecture r  and a remark:
   
    \begin{theo} (F. Zak, (Theorem 1, \cite{10}))
    
   Let  $X \subset  {\PP}^N$ be a   nondegenerate irreducible variety (not necessarily smooth). Consider a general projection $X \rightarrow {\PP}^{N-1}$. 
    If the quasi-projective locus  $X_k-X_{k+1}$  is connected and nonempty, then the hypersurfaces of degree $<k$ cut  complete linear system on $X$.    
    \end{theo}

 \begin{conj}   Let $X^n \subset {\PP}^N$ be a  dimension $n$, smooth irreducible variety. The following conditions are equivalent:
 
 1) the locus $X_k$ of a general projection of $X$ is reducible,
 
 2) the linear system cut out by the hypersurfaces of degree $k-1$ on $X$ is not complete.
 \end{conj}

 \begin{rem} The Palatini threefold (Example \ref{palvar})  is not quadratically normal.  The second author conjectured, at a Trento conference in $1988$,   that this is the only non quadratically normal smooth threefold in  $ {\PP}^5$.

 \end{rem}

\end{document}